\documentclass[a4paper,12pt]{article}
\usepackage{amsthm,amsmath}
\usepackage{amssymb}
\usepackage[latin1]{inputenc}
\usepackage[all]{xy}
\usepackage{enumerate}
\usepackage[T1]{fontenc}
\usepackage{a4wide}
\usepackage[mathscr]{eucal}
\usepackage{graphicx}
\usepackage{color}

\swapnumbers
\newtheorem{theorem}{Theorem}[section]
\newtheorem{proposition}[theorem]{Proposition}

\newtheorem{lemma}[theorem]{Lemma}
\newtheorem*{theorem*}{Theorem}
\newtheorem*{proposition*}{Proposition}
\newtheorem*{corollary*}{Corollary}
\newtheorem*{lemma*}{Lemma}
\theoremstyle{definition}
\newtheorem{definition}[theorem]{Definition}

\newtheorem{remark}[theorem]{Remark}
\newtheorem*{remark*}{Remark}

\newtheorem*{definition*}{Definition}

\newcommand{\cat}[1]{\mathcal{#1}}

\newcommand{\tensor}[1]{\otimes_{#1}}

\renewcommand{\hom}[3]{\mathrm{Hom}_{#1}(#2,\,#3)}

\newcommand{\bara}[1]{\overline{#1}}

\newcommand{\inv}[2]{\mathbf{Inv}_{#1}(#2) }

\newcommand{\fk}[1]{\mathfrak{#1}}
\newcommand{\Sf}[1]{\mathsf{#1}}

\newcommand{\lr}[1]{\left(\underset{}{} #1 \right)}

\newcommand{\Scr}[1]{\mathscr{#1}}
\newcommand{\Inv}[2]{\mathbf{Inv}_{#1}(#2)}
\newcommand{\Picar}[1]{\mathbf{Pic}(#1)}
\newcommand{\Aut}[2]{\mathbf{Aut}_{#1}(#2)}
\newcommand{\wid}[1]{\widehat{#1}}
\newcommand{\objeto}[3]{\xymatrix@C=35pt{ #1 \ar@2{->}|-{\,[#2]\,}[r]  & #3}}

\begin{document}
\allowdisplaybreaks

\title{Invertible unital bimodules over rings with local units, and related exact sequences of groups
\footnote{Research supported by  grant MTM2007-61673 from the Ministerio de
Educaci\'{o}n y Ciencia of Spain and FEDER, and P06-FQM-01889 from Junta de Andaluc\'{i}a}}
\author{L. El Kaoutit \\
\normalsize Departamento de \'{A}lgebra \\ \normalsize Facultad de Educaci\'{o}n y Humanidades \\
\normalsize  Universidad de Granada \\ \normalsize El Greco N${}^{0}$
10, E-51002 Ceuta, Espa\~{n}a \\ \normalsize
e-mail:\textsf{kaoutit@ugr.es} \and J. G\'omez-Torrecillas \\
\normalsize Departamento de \'{A}lgebra \\ \normalsize Facultad de Ciencias \\
\normalsize Universidad
de Granada\\ \normalsize E-18071 Granada, Espa\~{n}a \\
\normalsize e-mail: \textsf{gomezj@ugr.es} }

\date{\today}

\maketitle

\begin{abstract}
Given an extension $R \subseteq S$ of rings with same set of local units, inspired by the works of Miyashita, we construct four exact sequences of groups relating Picard's groups of $R$ and $S$.
\end{abstract}

\section*{Introduction}
Let $R \subseteq S$ be an extension of rings with identity element, and $\Scr{Z}(R)$, $\Scr{Z}(S)$ their centers. In   \cite{Miyashita:1970, Miyashita:1973}, Y\^{o}ichi Miyashita constructs the following four exact sequences of groups
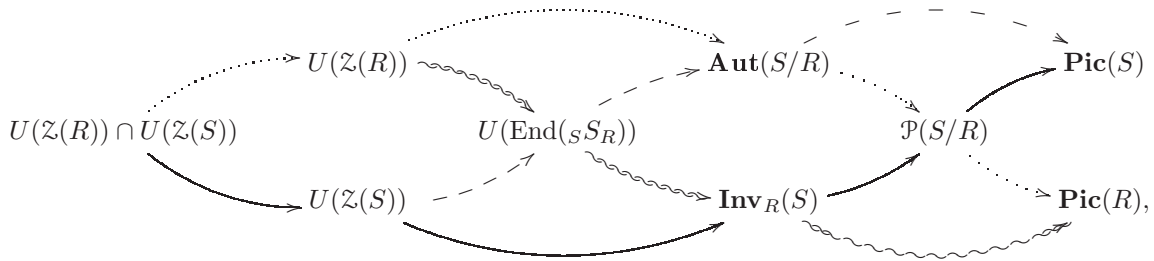
\begin{footnotesize}
$$
\xymatrixrowsep{10pt} \xymatrixcolsep{20pt}
\xymatrix{ & U(\Scr{Z}(R))  \ar@{.>}_-{}@/^1.7pc/[rr] \ar@{~>}_-{}@/^0.8pc/[rd] & & \textbf{Aut}(S/R) \ar@{.>}_-{}@/^0.6pc/[rd] \ar@{-->}_-{}@/^1.7pc/[rr] & &  \Picar{S} \\ 
U(\Scr{Z}(R)) \cap U(\Scr{Z}(S)) \ar@{..>}_-{}@/^1pc/[ru] \ar@{->}_-{}@/_1pc/[rd] & & U({\rm End}{({}_SS_R)}) \ar@{~>}_-{}@/_0.6pc/[rd] \ar@{-->}_-{}@/^0.6pc/[ru] & & \Scr{P}(S/R) \ar@{.>}_-{}@/_0.6pc/[rd] \ar@{->}_-{}@/^0.6pc/[ru] & \\ & U(\Scr{Z}(S)) \ar@{->}_-{}@/_1.7pc/[rr] \ar@{-->}_-{}@/_0.8pc/[ru] & & \Inv{R}{S} \ar@{->}_-{}@/_0.6pc/[ru] \ar@{~>}_-{}@/_1.7pc/[rr] & &  \Picar{R},}
$$
\end{footnotesize}
where $\textbf{Aut}(S/R)$ is the group of ring automorphisms of $S$ that acts by identity on $R$, $\Inv{R}{S}$ is the group of invertible $R$-subbimodules of $S$, $\Picar{R}$ and $\Picar{S}$ denote respectively the Picard group of $R$ and $S$; $\Scr{P}(S/R)$ is the group constructed in \cite[Theorem 1.3]{Miyashita:1973}, and wherein the notation $U(X)$ stands for the units group of an unital ring $X$. These sequences were the key steps used by Miyashita in his program of studying non-commutative Galois extensions using a generalized crossed product. One of the outcomes of that program was the generalization  to the non-commutative setting of the well known seven terms exact sequence given by Chase, Harrison and Rosenberg in \cite{Chase et al:1965}. He also proved that this seven terms exact sequence is a Morita-invariant object.

This paper is the first step of our project of extending Miyashita's program to the framework of small additive categories.  Given such a category $\cal{C}$, one can associate to it a ring with enough orthogonal idempotents $\Scr{R}(\cal{C}) =\oplus_{\fk{p},\, \fk{q} \, \in\, \cal{C}}\hom{\cal{C}}{\fk{p}}{\fk{q}}$ known as Gabriel's ring. If we think about all small additive categories and their bimodules as additive functors from $\cat{C}^{op}\times \cat{C}$ to the category of abelian groups, then the correspondence $\Scr{R}$ establishes in fact a bi-equivalence from the $2$-category of all small additive categories to the bicategory of unital bimodules over rings with enough orthogonal idempotents. This suggests that a most appropriate context for studying  small additive categories is that of unital bimodules over rings with enough orthogonal idempotents. Rings with local units (see Definition below) are  a slight generalization of rings with enough orthogonal idempotents. 
In this direction, our first task to be completed is the construction of sequences similar to those above for extension of rings $R \subseteq S$ with same set of local units. This will be fulfilled in this paper. It is worth noticing that Miyashita's contruction can not be transfered word by word to our context, since we are dealing with several objects rather than with only one.  This makes a big difference in time to make any calculations. It is noteworthy  that our methods can be seen as complete and detailed, even for the unital case.

\smallskip
{\textbf{Notations and basic notions}:}
In this paper ring means an associative and not necessary unital ring. We denote by $\Scr{Z}(R)$ (resp. $\Scr{Z}(G)$) the center of a ring $R$ (resp. of a group $G$), and if $R$ is unital we will denote by $U(R)$ the units group of $R$, that is, the set of all invertible elements of $R$.

Let $R$ be a ring with a fixed set of local units  $\Sf{E}$. This means that $\Sf{E}$  is a set of idempotent elements such that for every finite subset  $\{r_1,\cdots,r_n\}$ of $R$, there exists an element $ e \in\Sf{E}$ such that
$$
e\,r_i\,\,=\,\, r_i\, e\,\,=\,\, r_i, \quad \text{ for } \, i=1,\cdots,n,
$$
see \cite{Abrams:1983} and \cite{Anh/Marki:1983}. We will use the notation $$Unit\{r_1,\cdots,r_n\}\,\,:=\,\, \left\{\underset{}{} e \in\Sf{E}|\,\, er_i= r_ie= r_i, \,\, \text{for} \, i=1,2,\cdots,n  \right\}.$$

A right $R$-module $X$ is said to be \emph{unital} provided one of the following equivalent conditions
\begin{enumerate}[(i)]
\item $X\tensor{R}R \cong X$ via the right $R$-action of $X$,
\item $XR:=\{\sum_{finite} x_ir_i|\, r_i \in R,\, x_i \in X\}=X$,
\item for every element $x \in X$, there exists an element $e \in \Sf{E}$ such that $xe=x$.
\end{enumerate}

Left unital $R$-modules are  defined analogously. A unital $R$-bimodule is an $R$-bimodule which is left and right unital. Obviously any right $R$-module $X$ contains $XR$ as the largest right unital $R$-submodule.
In all what follows we consider an extension of rings with local units $R \subseteq S$ (having the same set of local unit $\Sf{E}$). Observe that since $R$ and $S$ have the same set of local units, any right (resp. left) unital $S$-module can be considered as right (resp. left) unital $R$-module by restricting scalars. Furthermore, for any right $S$-module $X$, we have $XR=XS$.
Given two subsets $W$ and $V$  of $S$, we denote by
$$
WV\,\,=\,\, \left\{\underset{}{}  \sum_{\text{finite}} w_iv_i \in S|\,\, w_i \in W,\, v_i \in V  \right\}.
$$

\section{The first exact sequence of groups}\label{Sec:1}

Let $R \subseteq S$ be an extension of  rings with local units, we denote by $\Inv{R}{S}$ the set of invertible unital $R$-subbimodules of $S$. That is, $X \in \Inv{R}{S}$ provided that $X$ is a unital $R$-subbimodule of $S$ and there exists another $R$-subbimodule $Y$ of $S$ such that $ XY=YX= R$. Clearly, this defines a group structure on $\inv{R}{S}$ whose neutral element is $R$.

\begin{lemma}\label{lema:1}
Let $X$ be an element of the group $\inv{R}{S}$ with inverse $Y$. Then there are two $R$-bilinear isomorphisms
$$
\xymatrix{ X\tensor{R}Y \ar@{->}^-{\cong}[r] & R & Y \tensor{R}X. \ar@{->}_-{\cong}[l]}
$$ In particular, the maps
$$
\xymatrix@R=0pt{S\tensor{R}X \ar@{->}^-{m_l}[r] & S, \\ s\tensor{R}x \ar@{|->}[r] & sx } \text{ and } \,\, \xymatrix@R=0pt{ X\tensor{R}S \ar@{->}^-{m_r}[r] & S \\ x\tensor{R}s \ar@{|->}[r]  &  xs }
$$
are isomorphisms of $S$-$R$ and $R$-$S$-bimodules, respectively.
\end{lemma}
\begin{proof}
The first isomorphism is given by the map $\varsigma: R \longrightarrow X\tensor{R}Y$ which sends $ r\longmapsto \sum_{(r)} x_{(r)}\tensor{R}y_{(r)}$, where $r\,=\, \sum_{(r)}x_{(r)}y_{(r)} \in R=XY$. This is a well defined map since $\sum_{(r)} x_{(r)}\tensor{R}y_{(r)} \,=\,0$ if and only if $r=\sum_{(r)}x_{(r)}y_{(r)}\,=\, 0$. This equivalence says that $\varsigma$ is also injective.  Now, given an element $r \in R$, we can easily check using the equalities $R\,=\, XY\,=\, YX$, that $\varsigma(r)\,=\,\varsigma(e)r\,=\,r\varsigma(e)$, for any unit $e \in Unit\{r\}$ and every element $r \in R$. Thus $\varsigma$ is an $R$-bilinear isomorphism. Similarly we show that $R \cong Y\tensor{R}X$. The inverse map of $m_r$ is defined as follows: Given an element $s  \in S$, with unit $e$, we set $m_r^{-1}(s)\,=\, \sum_{(e)}x_{(e)}\tensor{R}y_{(e)}s \, \in X\tensor{R}S$, where $e \,=\, \sum_{(e)}x_{(e)}y_{(e)} \in XY$. If $f$ is another unit for $s$, then we can consider $h \in Unit\{e,f\}$ and we have
\begin{eqnarray*}
  \sum_{(h)}x_{(h)}\tensor{R}y_{(h)}s &=& \sum_{(h)}x_{(h)}\tensor{R}y_{(h)}es  \\
   &=& \sum_{(h),\, (e)}x_{(h)}\tensor{R}y_{(h)}x_{(e)}y_{(e)}s \\
   &=& \sum_{(h),\, (e)}x_{(h)}y_{(h)}x_{(e)}\tensor{R}y_{(e)}s \\
   &=& \sum_{ (e)}h\,x_{(e)}\tensor{R}y_{(e)}s,\quad \sum_{(e)}x_{(e)}\tensor{R}y_{(e)}h\,=\,\sum_{(e)}hx_{(e)}\tensor{R}y_{(e)} \\
   &=& \sum_{(e)}x_{(e)}\tensor{R}y_{(e)}\,h\,s, \quad s=hs=sh \\
   &=& \sum_{(e)}x_{(e)}\tensor{R}y_{(e)}s.
\end{eqnarray*}
By the same way, we obtain the equality $\sum_{(h)}x_{(h)}\tensor{R}y_{(h)}s\,=\, \sum_{(f)}x_{(f)}\tensor{R}y_{(f)}s$. Whence $\sum_{(e)}x_{(e)}\tensor{R}y_{(e)}s\,=\, \sum_{(f)}x_{(f)}\tensor{R}y_{(f)}s$, and $m_r^{-1}$ is a well defined inverse map of $m_r$. Similarly, we show that $m_l$ is an isomorphism of unital $(S,R)$-bimodules.
\end{proof}

The Picard group $\Picar{R}$ of a ring with local units $R$ is defined to be the set of isomorphism classes $[P]$ where $P$ is a unital $R$-bimodule such that there exists a unital $R$-bimodule $Q$ with $R$-bilinear isomorphisms 
\begin{equation}\label{Picard-iso}
\xymatrix{Q\tensor{R}P \ar@{->}^-{\fk{r}}_-{\cong}[r] & R &  P\tensor{R}Q. \ar@{->}_-{\fk{l}}^-{\cong}[l]}
\end{equation}
These isomorphisms define in fact an auto-equivalence on the category of unital right $R$-modules via the tensor product functors $-\tensor{R}P$ and $-\tensor{R}Q$. They also satisfy the following equalities
\begin{eqnarray}
  (P\tensor{R}\fk{r}) \circ (P\tensor{R}Q\tensor{R}\fk{l}\tensor{R}P) &=& (\fk{l}\tensor{R}P) \circ (P\tensor{R}\fk{r}\tensor{R}Q\tensor{R}P) \label{lr-1} \\
  (\fk{r}\tensor{R}Q) \circ (Q\tensor{R}\fk{l}\tensor{R}P\tensor{R}Q) &=&  (Q\tensor{R}\fk{l}) \circ (\fk{r}\tensor{R}Q\tensor{R}P\tensor{R}Q). \label{lr-2}
\end{eqnarray}
Furtheremore, we can deduce from \cite[Proposition 5.1]{Kaoutit:2006} that the right (resp. left) unital $R$-module $eP$ (resp. $Pe$) is finitely generated and projective, for every unit $e \in \Sf{E}$. 
Thus, for every right unital $R$-module $N$,  there is a natural isomorphism $$\Theta_N: N\tensor{R}Q \longrightarrow \hom{R}{P}{N}R$$ defined as follows. For every $n \in N$, $q \in Q$, and $p \in P$ with a common unit $e \in \Sf{E}$, we have
\begin{equation}\label{Theta}
\Theta_N(n\tensor{R}q)(p)\,=\,(N\tensor{R}\fk{r}) \circ (N\tensor{R}Q\tensor{R}\fk{l}\tensor{R}P) \circ (N\tensor{R}\fk{r}^{-1}\tensor{R}Q\tensor{R}P) (n\tensor{R}e \tensor{R}q\tensor{R}p).
\end{equation}
\begin{remark}
An element $[P] \in \Picar{R}$ is not necessarily represented by a finitely generated right (or left) module $P_R$ (or ${}_RP$). It is in fact a direct limit of finitely generated and projective right modules of the form $Pe$, $e \in \Sf{E}$. This of course marks a big difference between the study of the Picard group of a ring with unit and a ring with local units.
\end{remark}

The multiplication in $\Picar{R}$ is defined by $[P] [P'] = [P\tensor{R}P']$ and the neutral element is $[R]$. By Lemma \ref{lema:1}, we have a morphism of groups
$$ \xymatrix{\inv{R}{S} \ar@{->}^-{[-]}[rr] & & \Picar{R}.}$$
Now, define $\mathbf{Aut}_{S-R}(S)$ to be the units group of the monoid $\hom{S-R}{S}{S}$ of all $S$-$R$-bilinear endomorphisms of $S$. In other words the set of all $S$-$R$-bilinear automorphisms of $S$. In the unital case one can easily check that this group is isomorphic to the group of units of the subring of $R$-invariant elements of $S$. Here, one can show that for every idempotent element $e \in R $ and every element $\lambda \in \Aut{S-R}{S}$, we have
\begin{equation}\label{lambda-E}
\lambda(e) \lambda^{-1}(e) \,\,=\,\, e \,\, =\,\, \lambda^{-1}(e) \lambda(e).
\end{equation}
That is, $\lambda(e) $ belongs to the units group of the unital subring of $(eRe)$-invariant elements of $eSe$.

\begin{lemma}\label{lema:2}
Let $R \subseteq S$ be a ring extension with the same set of local units. Then the map
$$\xymatrix@R=0pt{ \Aut{S-R}{S} \ar@{->}^-{\Scr{D}}[rr] & & \inv{R}{S} \\ \lambda \ar@{|->}[rr]  & & \lambda^{-1}(R) }$$ is a homomorphism of groups with  kernel in the center sub-group of $\Aut{S-R}{R}$, i.e.,  $${\rm Ker}(\Scr{D}) \subseteq \Scr{Z}(\Aut{S-R}{R}).$$
\end{lemma}
\begin{proof}
For every $\lambda \in \Aut{S-R}{R}$, we have
$$\lambda(R)\lambda^{-1}(R)\,=\, \lambda^{-1}(\lambda(R)R) \,=\, \lambda^{-1} \circ \lambda(R)\,=\,R.$$
Similarly, we get $\lambda^{-1}(R)\lambda(R)=R$. Therefore, $\lambda^{-1}(R) \in \inv{R}{S}$. Thus $\Scr{D}$ is a well defined map. Now, we have  $\Scr{D}(1)=R$ (here $1$ denotes the neutral element of the group $\Aut{S-R}{R}$), and for every pair of elements $\lambda$ and $\gamma$ in $\Aut{S-R}{S}$, we have
\begin{eqnarray*}
 \Scr{D}(\lambda \circ \gamma)  &=& \gamma^{-1}(\lambda^{-1}(R)) \\
   &=& \gamma^{-1}(\lambda^{-1}(R)R) \\
   &=& \lambda^{-1}(R)\gamma^{-1}(R) \\
   &=& \Scr{D}(\lambda) \Scr{D}(\gamma).
\end{eqnarray*}
Thus $\Scr{D}$ is a morphism of groups. Moreover, if $\lambda \in {\rm Ker}(\Scr{D})$, then $\lambda(R)=R$. So let $\gamma$ be an arbitrary element of $\Aut{S-R}{S}$, we have \begin{eqnarray*}
 \lambda \circ \gamma(s)  &=& \lambda (\gamma(se)), \quad es=se=s, \,\, e \in Unit\{s\} \\
   &=&  \gamma(s)\lambda(e),\quad \lambda(e) \in R\\
   &=&  \gamma(s\lambda(e)) \\
   &=& \gamma \lambda (se) \,\, =\,\, \gamma \circ \lambda (s),
\end{eqnarray*}
for every element $s \in S$. So we have $\lambda \circ \gamma \,\,=\,\, \gamma \circ \lambda$, for every element $ \gamma \in \Aut{S-R}{S}$. Whence, ${\rm Ker}(\Scr{D}) \subseteq \Scr{Z}(\Aut{S-R}{S})$.
\end{proof}

\begin{proposition}
Let $R \subseteq S$ be a ring extension with the same set of local units $\Sf{E}$. Then we have a commutative diagram of groups with exact rows
$$
\xymatrix{ 1 \ar@{->}^-{}[r] & {\rm Ker}(\Scr{D}) \ar@{_{(}->}^-{}[d] \ar@{->}^-{}[r] & \Aut{S-R}{R} \ar@{=}^-{}[d] \ar@{->}^-{\Scr{D}}[r] & \inv{R}{S} \ar@{->}^-{[-]}[r] & \Picar{R}  \\ 1 \ar@{->}^-{}[r] & \Scr{Z}(\Aut{S-R}{R}) \ar@{->}^-{}[r] &  \Aut{S-R}{R} & &  }
$$
Moreover, we have
$$
{\rm Ker}(\Scr{D})\,\,=\,\, \left\{\underset{}{}  \lambda \in \Aut{S-R}{S}|\,\, \lambda(e) \in U(\Scr{Z}(eRe)), \text{ for every element } e \in \Sf{E} \right\}.
$$
\end{proposition}
\begin{proof}
Let us first show that ${\rm Ker}([-]) = {\rm Im}(\Scr{D})$. So given $X \in \inv{R}{S}$ such that $[X]=[R]$, we know that there is an $R$-bilinear isomorphism $\varphi: X \to R$. Set
$$\xymatrix{ \lambda: S \ar@{->}^-{m_l^{-1}}[rr] & & S\tensor{R}X \ar@{->}^-{S\tensor{R}\varphi}[rr] & & S\tensor{R}R \cong S, }$$
where $m_l$ is the isomorphism of Lemma \ref{lema:1}. By definition $\lambda$ is an  $S$-$R$-bilinear automorphism of $S$. We claim that $\lambda^{-1}(R)\,=\,X$. To this end, consider $r \in R$ and put $r\,=\,\sum_{(r)}y_{(r)}x_{(r)} \in YX$, where $Y$ is the inverse of $X$ in $\inv{R}{S}$. So $\lambda(r)\,=\, \sum_{(r)}y_{(r)} \varphi(x_{(r)})$, and we have
\begin{eqnarray*}
 \lambda^{-1}(r)  &=& m_l \circ (S\tensor{R}\varphi^{-1}) (r\tensor{R}e),\quad \text{ where }\, r=er=re,\,e \in \Sf{E}  \\
   &=& m_l(r\tensor{R}\varphi^{-1}(e)) \\
   &=&  r\varphi^{-1}(e) \in X,
\end{eqnarray*}
and this shows that $\lambda^{-1}(R) \subseteq X$. Conversely, given an element $x \in X$, we have $\varphi(x) \in R$, and
\begin{eqnarray*}
    \lambda^{-1}(\varphi(x))   &=& m_l \circ (\varphi(x)\tensor{R}\varphi^{-1}(e)),\quad \text{ where }\, ex=xe=x,\,e \in \Sf{E}  \\
       &=& \varphi(x)\varphi^{-1}(e) \\
       &=&  \varphi^{-1}(\varphi(x)e)\,\,=\,\, \varphi^{-1} \circ \varphi(x)\,\,=\,\, x,
\end{eqnarray*}
that is, $X \subseteq \lambda^{-1}(R)$ which completes the proof of the claim. This shows that ${\rm Ker}([-]) \subseteq {\rm Im}(\Scr{D})$. Conversely, let $X \in { \rm Im}(\Scr{D})$, that is, $X= \lambda^{-1}(R)$ for some $\lambda \in \Aut{S-R}{S}$. The restriction of $\lambda$ to the $R$-subbimodule $R$ gives an $R$-bilinear isomorphism $\varphi=\lambda^{-1}_{/R}: R \to X$. This means that $[X]=[R]$ in $\Picar{R}$. Therefore, ${\rm Im}(\Scr{D}) \subseteq {\rm Ker}([-])$, and this proves the stated exactness. The first statement is then deduced from Lemma \ref{lema:2}.

Now, let $\lambda \in \Aut{S-R}{S}$ such that $\Scr{D}(\lambda)=R$, that is, $\lambda^{-1}(R)=R$. So we have $\lambda(R)=R$, and that $\lambda(e) \in R$ for every unit $e \in \Sf{E}$. Given an arbitrary element $r \in R$ and $e \in Unit\{r\}$, we have $\lambda(r)=r\lambda(e)=\lambda(e)r$. Hence, $\lambda(e) \in \Scr{Z}(eRe)$, for every element $e \in \Sf{E}$. By equation \eqref{lambda-E}, we get that $\lambda(e) \in U(\Scr{Z}(eRe))$, for every element $e \in \Sf{E}$. Conversely, if $\lambda \in \Aut{S-R}{S}$ such that $\lambda(e) \in U(\Scr{Z}(eRe))$, for every  element  $ e \in \Sf{E}$, then clearly we have $\lambda(R) \subseteq  R$. The reciprocal inclusion also holds since by equation \eqref{lambda-E}, $\lambda^{-1}(r) = r \lambda^{-1}(e) \in R$, for every $r \in R$ and $e \in Unit\{r\}$. Thus $\lambda(R)=R$, and the stated equality is proved.
\end{proof}

\section{The second exact sequence of groups}\label{Sect-II}

Let us  denote by $\Aut{R-rings}{S}$ the units group of the monoid ${\bf End}_{R-rings}(S)$ of all ring endomorphisms which act by identity on  $R$.  Given an element $\phi \in \Aut{R-rings}{S}$ and unital $S$-bimodule $M$, we denote by $M_{\phi}$ the unital $S$-bimodule whose underlying
left $S$-module coincides with ${}_SM$ and its right $S$-module structure is given by restriction of scalars using $\phi$. That is, $M$ with the following $S$-biaction
$$ t. m_{\phi}  \,\,=\,\, (tm)_{\phi},\,\,\text{ and }\,\, m_{\phi}.t\,\,=\,\, (m\phi(t))_{\phi}, \,\,\text{for every }\,\, m_{\phi} \in M_{\phi},\, t \in S.$$

\begin{lemma}\label{lema:21}
Let $R \subseteq S$ be a ring extension with the same set of local units $\Sf{E}$. Then there is a morphism of groups
$$
\xymatrix@R=0pt{ \Aut{R-rings}{S} \ar@{->}^-{[S_{-}]}[rr] & & \Picar{S} \\ \phi \ar@{|->}[rr] & & S_{\phi}.   }
$$
\end{lemma}
\begin{proof}
It is clear that $[S_{1}]=[S]$, where $1$ is the neutral element of the group $\Aut{R-rings}{S}$. Now, if $\phi, \psi \in \Aut{R-rings}{S}$, then the identity map $S_{\phi\psi} \to (S_{\phi})_{\psi}$ sending $s_{\phi\psi} \mapsto (s_{\phi})_{\psi}$ defines an $S$-bilinear isomorphism. Namely, the left $S$-action stills by definition inalterable, so we need to check the right one. Given an element $t \in S$, we have
$$ s_{\phi\psi} . t\,\,=\,\, s\phi\psi(t)\,\,=\,\, s_{\phi} . \psi(t)\,\,=\,\, (s_{\phi})_{\psi} . t,$$ for every element $s_{\phi\psi} \in S_{\phi\psi}$. This shows that $[S_{\phi\psi}]=[(S_{\phi})_{\psi}]$ in the group $\Picar{S}$. On the other hand, we know that $S_{\phi}$ is a right unital $S$-module, thus $S_{\phi} \cong S_{\phi}\tensor{S}S$. Applying the restricting scalars functor $(-)_{\psi}$, we obtain $(S_{\phi})_{\psi} \cong (S_{\phi}\tensor{S}S)_{\psi}$. We then get $[S_{\phi\psi}]=[S_{\phi}\tensor{S}S_{\psi}]$, since the map $(S_{\phi}\tensor{S}S)_{\psi} \to S_{\phi}\tensor{S}S_{\psi}$ sending $(s_{\phi}\tensor{S}t)_{\psi} \mapsto s_{\phi}\tensor{S}t_{\psi}$ is an $S$-bilinear isomorphism. We then conclude that the stated map is a morphism of groups.
\end{proof}

\begin{lemma}\label{lema:22}
Let $R \subseteq S$ be a ring extension with the same set of local units $\Sf{E}$. Then there is a morphism of groups
$$
\xymatrix@R=0pt{  \Aut{S-R}{S} \ar@{->}^-{\wid{(-)}}[rr] & & \Aut{R-rings}{S} \\ \lambda \ar@{|->}[rr] & &  \wid{\lambda}, }
$$ where $\wid{\lambda}: S \longrightarrow S$ sends $s \longmapsto \lambda^{-1}(e) \,s\, \lambda(e)$, for $e \in Unit\{s\}$.
\end{lemma}
\begin{proof}
Let us first check that $\wid{\lambda}$ is a well defined map. Consider $f,e \in Unit\{s\}$ two units for a fixed element $s$, and let $h \in Unit\{e,f\}$. Then, we have
\begin{eqnarray*}
 \lambda^{-1}(h) s \lambda(h)  &=& \lambda^{-1}(h) ese\lambda(h)  \\
   &=& \lambda^{-1}(he) s \lambda(eh)  \\
   &=&  \lambda^{-1}(e) s \lambda(e).
\end{eqnarray*}
Similarly, we get $\lambda^{-1}(h) s \lambda(h)=\lambda^{-1}(f) s \lambda(f)$. Thus $\wid{\lambda}$ is independent from the chosen unit. This map is clearly additive, and for every pair of elements  $s,t \in S$, we have
\begin{eqnarray*}
 \wid{\lambda}(st)  &=& \lambda^{-1}(e)\, s\, t\, \lambda(e),\quad e \in Unit\{s,t\}   \\
 &=& \lambda^{-1}(e)\, s\,e\, t\, \lambda(e)  \\
   &\overset{\eqref{lambda-E}}{=}& \lambda^{-1}(e)\, s\,\lambda(e) \lambda^{-1}(e) t\, \lambda(e) \\
   &=&  \wid{\lambda}(s) \wid{\lambda}(t).
\end{eqnarray*}
This shows that $\wid{\lambda}$ is multiplicative. Now, given an element $r \in R$, we have
\begin{eqnarray*}
 \wid{\lambda}(r)  &=& \lambda^{-1}(e) r \lambda(e),\quad e \in Unit\{r\} \\
   &=& \lambda^{-1}(er)\lambda(e) \\
   &=&  \lambda^{-1}(r)\lambda(e) \\
   &=& \lambda\lambda^{-1}(re)\,\,=\,\, \lambda\lambda^{-1}(r)\,\,=\,\,r.
\end{eqnarray*}
We have proved that $\wid{\lambda} \in \Aut{R-rings}{S}$. It is clear that $\wid{1}=1$. Given a pair of elements $\gamma, \lambda \in \Aut{S-R}{S}$, we have
\begin{eqnarray*}
 \wid{\lambda\circ \gamma}(s)  &=& (\lambda\circ \gamma)^{-1}(e) \, s\, (\lambda \circ \gamma)(e),\quad e \in Unit\{s\} \\ &=& \gamma^{-1}(\lambda^{-1}(e)) \, s \, \lambda(\gamma(e)) \\
   &=& \lambda^{-1}(e)\gamma^{-1}(e) \,s \, \gamma(e)\lambda(e) \\
   &=& \lambda^{-1}(e)\wid{\gamma}(s) \lambda(e) \\
   &=&  \wid{\lambda} \circ \wid{\gamma}(s).
\end{eqnarray*}
Therefore, $\wid{(-)}$ is a morphism of groups.
\end{proof}

\begin{proposition}\label{pro:21}
Let $R \subseteq S$ be a ring extension with the same set of local units $\Sf{E}$. Then there is an exact sequence of groups
$$
\xymatrix@C=40pt{ 1 \ar@{->}^-{}[r] & {\rm Ker}\lr{\wid{(-)}} \ar@{->}^-{}[r] & \Aut{S-R}{S} \ar@{->}^-{\wid{(-)}}[r] & \Aut{R-rings}{S} \ar@{->}^-{[S_{-}]}[r] & \Picar{S},  }
$$
where $${\rm Ker}\lr{\wid{(-)}}\,\,=\,\,\left\{\underset{}{} \lambda \in \Aut{S-R}{S}|\, \lambda(e) \in U(\Scr{Z}(eSe)), \text{ for every idempotent }\, e \in \Sf{E}  \right\}.$$ In particular ${\rm Ker}\lr{\wid{(-)}}\,=\,\Aut{S-S}{S}$ the group of all $S$-bilinear automorphisms of $S$.
\end{proposition}
\begin{proof}
Let $\phi \in \Aut{R-rings}{S}$ such that $[S_{\phi}] = [S]$ in the Picard group $\Picar{S}$. This implies that there exists an $S$-linear isomorphism $\omega: {}_SS_S \to  {}_S(S_{\phi})_S $. In particular, this map satisfies the following equation
\begin{equation}\label{omega}
\omega(t\,s\,u)\,\,=\,\, t\,\omega(s)\, \phi(u),\quad\text{ for every elements }\,\, s,t, u \in S.
\end{equation}
Since $\phi$ acts by identity on $R$, the isomorphism $\omega$ induces an $S$-$R$-bilinear isomorphism $\lambda: {}_SS_R \to {}_SS_R$ sending $s \mapsto \omega(s)$. Given $s \in S$ with unit $e \in Unit\{s\}$, we have
\begin{eqnarray*}
  \wid{\lambda}(s) &=& \lambda^{-1}(e) \,s\, \lambda(e)  \\
   &=& \omega^{-1}(e)\, s\, \omega(e) \\
   &\overset{\eqref{omega}}{=}& \omega^{-1}(e\phi(s))\, \omega(e) \\
   &=&  \omega\lr{\omega^{-1}(e\phi(s))\,e} \\
   &=& \omega\omega^{-1}\lr{\phi(e)\phi(s)\phi(e)}\,\,=\,\, \phi(ese)\,\,=\,\, \phi(s),
\end{eqnarray*}
that is, $\wid{\lambda}\,=\, \phi$, and this shows that ${\rm Ker}([S_{-}]) \subseteq {\rm Im}(\wid{(-)})$.

Conversely, given $\lambda \in \Aut{S-R}{S}$, consider the unital $S$-bimodule ${}_S(S_{\wid{\lambda}})_S$. Define a map $\omega: {}_S(S_{\wid{\lambda}})_S \to {}_SS_S$ by sending $s \mapsto \lambda^{-1}(s)$. This is clearly a left $S$-linear map. On the other hand, if $s, t \in S$ with $e \in Unit\{s,t\}$ are given, then we have
\begin{eqnarray*}
 \omega(s\, \wid{\lambda}(t))  &=& \lambda^{-1}\lr{s\, \lambda^{-1}(e)\,t\,\lambda(e)}  \\
   &=& \lambda^{-1}\lr{\lambda^{-1}(s)\,\lambda(t)} \\
   &=& \lambda^{-1}(s) \, \lambda^{-1}\lr{\lambda(t)}\\
   &=& \lambda^{-1}(s)\, t\,\,=\,\, \omega(s) t,
\end{eqnarray*}
which implies that $\omega$ is also right $S$-linear. Therefore, $[S_{\wid{\lambda}}]=[S]$ in $\Picar{S}$, for every $\lambda \in \Aut{S-R}{S}$, and this proves that the stated sequence in exact.

Now, let $\lambda \in {\rm Ker}(\wid{(-)})$. Then for every  $s \in S$, we have $\wid{\lambda}(s)=s$. By equation \eqref{lambda-E}, we thus get $\lambda(e) s=s \lambda(e)$, for every element $s \in S$ with $e \in Unit\{s\}$. In particular, we have $$ ete \lambda(e) \,=\, \lambda(e) ete,\,\,  \text{ for every } e \in \Sf{E} \text{ and } t \in S.$$ Therefore, $\lambda(e) \in U(\Scr{Z}(eSe))$, for every  $e \in \Sf{E}$. Conversely, consider  $\lambda \in \Aut{S-R}{S}$ such that $\lambda(e) \in U(\Scr{Z}(eSe))$, for every  $e \in \Sf{E}$. Then, we have $$\wid{\lambda}(s)=\lambda^{-1}(e)\,s\, \lambda(e) \,\,=\,\, \lambda^{-1}(e)\lambda(e)\, s\,\,=\,\, e \,s\,\,=\,\,s,$$ where the third equality is derived from equation \eqref{lambda-E}, and thus
$${\rm Ker}(\wid{(-)}) \,\,=\,\, \left\{\underset{}{} \lambda \in \Aut{S-R}{S}|\, \lambda(e) \in U(\Scr{Z}(eSe)), \text{ for every idempotent }\, e \in \Sf{E}  \right\}.$$ The last stated equality is an easy consequence of the previous equality.
\end{proof}

\section{The group $\Scr{P}(S/R)$}\label{Sect-III}
Let $R \subseteq S$ be an extension of rings with the same set of local units $\Sf{E}$. Given a unital $S$-bimodule $X$, we can consider it obviously as a unital $R$-bimodule. Let $P$ is a unital $R$-bimodule and $X$ a unital $S$-bimodule together with an $R$-linear map $\phi: P \to X$ such that
\begin{equation}\label{Phi-l-r}
\xymatrix@R=0pt{P\tensor{R}S \ar@{->}^-{\bara{\phi}_r}[r] & X, \\ p\tensor{R}s \ar@{|->}[r] & \phi(p)s} \qquad \xymatrix@R=0pt{S\tensor{R}P \ar@{->}^-{\bara{\phi}_l}[r] & X \\ t\tensor{R}q \ar@{|->}[r] & t\phi(q) } 
\end{equation}
are isomorphisms, respectively, of $(R,S)$-bimodules and $(S,R)$-bimodules. We will denote this situation by $\objeto{P}{\phi}{X}$. These objects define in fact a category which we denote by $\Scr{M}(S/R)$ and given by the following data:
\begin{description}
  \item[{\emph{Objects:}}]They are those $\objeto{P}{\phi}{X}$ described above.
  \item[\emph{Morphisms:}]They are pairs $(\alpha,\beta): (\objeto{P}{\phi}{X}) \,\to\, (\objeto{P'}{\phi'}{X'})$ where $\alpha: P \to P'$ and $\beta: X \to X'$ are, respectively, $R$-bilinear  and $S$-bilinear map, rendering commutative the following diagram
      $$\xy  *+{P}="p", p+<4cm, 0pt>*+{X}="1", p+<0pt,-1.5cm>*+{P'}="2", p+<4cm,-1.5cm>*+{X'}="3", {"p" \ar@{->}^-{\phi} "1"}, {"p" \ar@{->}_-{\alpha} "2"}, { "1" \ar@{->}^-{\beta} "3"}, {"2" \ar@{->}_-{\phi'}"3"} \endxy$$
\end{description}

The composition operation and the identity morphisms are the obvious ones. This category is in fact a monoidal category with multiplication defined as follows: Given two objects $\objeto{P}{\phi}{X}$ and $\objeto{Q}{\psi}{Y}$ in $\Scr{M}(S/R)$, one can define the map $$\xymatrix@C=40pt{\chi: P\tensor{R}Q \ar@{->}^-{\phi\tensor{R}\psi}[r] & X\tensor{R}Y \ar@{->}^-{\omega_{X,\,Y}}[r] & X\tensor{S}Y,}$$ where $\omega_{-,-}$ is the obvious natural transformation. The map $\chi$ is clearly $R$-bilinear. By definition, we have a chain of isomorphisms
\begin{footnotesize}
$$
\xymatrix@C=45pt{ P\tensor{R}Q\tensor{R}S \ar@{-->}_-{\bara{\chi}_r}[rrdd] \ar@{->}^-{P\tensor{R}\bara{\psi}_r}[r] & P\tensor{R}Y \ar@{->}^-{P\tensor{R}\bara{\psi}_l{}^{-1}}[r] & P\tensor{R}S\tensor{R}Q \ar@{->}^-{\cong}[d]  \\  & &  P\tensor{R}S\tensor{S}S\tensor{R}Q \ar@{->}^-{\bara{\phi}_r\tensor{S}\bara{\psi}_l }[d] \\ & & X\tensor{S}Y,  }
$$
\end{footnotesize}
whose composition is exactly $\bara{\chi}_r$. Similarly we show that $\bara{\chi}_l$ is also an isomorphism. This proves that $\objeto{(P\tensor{R}Q)}{\chi}{(X\tensor{S}Y)}$ is again an object of the category $\Scr{M}(S/R)$ which we denote  by $\xymatrix@C=50pt{ (P\tensor{R}Q) \ar@{=>}|-{\,[\phi] \cdot[\psi]\,}[r]  & (X\tensor{S}Y)}$. 

Now, given two morphisms  $(\alpha,\beta): (\objeto{P}{\phi}{X}) \to (\objeto{P'}{\phi'}{X'})$ and $(\mu,\nu): (\objeto{Q}{\psi}{Y}) \to (\objeto{Q'}{\psi'}{Y'})$ in $\Scr{M}(S/R)$, we have a morphism  $$\xymatrix@R=30pt@C=40pt{ (P\tensor{R}Q) \ar@2{->}|-{[\phi]\cdot[\psi]}[rr] & \ar@{->}|-{(\alpha\tensor{R}\mu,\,\beta\tensor{R}\nu)}[d] & (X\tensor{S}Y) \\ (P'\tensor{R}Q') \ar@2{->}|-{[\phi']\cdot[\psi']}[rr] & & (X'\tensor{S}Y')  }$$ in $\Scr{M}(S/R)$. The unit object with respect to this multiplication is proportioned by the object $\objeto{R}{\iota}{S}$, where $\iota$ is the inclusion $R \subseteq S$. The axioms of a monoidal category are easily verified for $\Scr{M}(S/R)$ since it was built upon the monoidal categories of unital $R$-bimodules and unital $S$-bimodules.

\begin{lemma}\label{lema:31}
Let $R$ and $S$ be as above, and consider a unital $R$-bimodule $P$ and a unital $S$-bimodule $X$ with an $R$-bilinear map $\phi: P \to X$ with associated maps $\bara{\phi}_r$ and $\bara{\phi}_l$ as in \eqref{Phi-l-r}. If $[P] \in \Picar{R}$ and $[X] \in \Picar{S}$, then the following are equivalent
\begin{enumerate}[(i)]
\item $\bara{\phi}_r: P\tensor{R}S \to X$ is an $R$-$S$-bilinear isomorphism;
\item $\bara{\phi}_l: S\tensor{R}P \to X$ is an $S$-$R$-bilinear isomorphism.
\end{enumerate}
Moreover, $\phi$ is an injective map provided one of the conditions $(i)$ or $(ii)$.
\end{lemma}
\begin{proof}
We only prove $(i) \Rightarrow (ii)$, the reciprocal implication follows similarly. So assume that $\bara{\phi}_r$ is an isomorphism, since  $[X] \in \Picar{S}$, we have the following chain of isomorphisms:
\begin{footnotesize}
$$
\xymatrix@C=30pt{S \ar@{-->}_-{\kappa}[rrrdd] \ar@{->}^-{\cong}[r] & S\,\hom{S}{X}{X} \ar@{->}^-{\cong}[r] & S\,\hom{S}{P\tensor{R}S}{X} \ar@{->}^-{\cong}[r]
&  S\,\hom{R}{P}{\hom{R}{S}{X}} \ar@{=}|-{}[d] \\ & & & S\,\hom{R}{P}{\hom{R}{S}{X}R} \ar@{->}|-{\cong}[d] \\ & & & S\,\hom{R}{P}{X}, }
$$
\end{footnotesize}
whose composition $\kappa$ is given explicitly by $\kappa(s): P \to X$ sending $p \mapsto s\,\phi(p)$, for every element $s \in S$. Now, since $[P] \in \Picar{R}$, we know that $eP$ is a finitely generated and projective right unital $R$-module, for every $e \in \Sf{E}$. Therefore, the map
$$  
\upsilon: S\, \hom{R}{P}{X} \tensor{R}P \longrightarrow X,\quad \lr{sf\tensor{R}p \longmapsto  sf(p)}
$$
is an $S$-$R$-bilinear isomorphism. This implies that $\bara{\phi}_l = \upsilon \circ (\kappa\tensor{R}P)$ is also an $S$-$R$-bilinear isomorphism. The rest of the proof is clear since $P_R$ and ${}_RP$ are flat modules.
\end{proof}

\begin{definition}\label{PSR}
Let $R \subseteq S$ be an extension of rings with the same set of local units, and consider its associated monoidal category $\Scr{M}(S/R)$. Let us denote by $\Scr{S}(S/R)$ the skeleton of $\Scr{M}(S/R)$ (i.e., the set of all isomorphism classes). We consider $\Scr{S}(S/R)$ as a monoid with multiplication and unit induced by the monoidal structure of $\Scr{M}(S/R)$. We define $\Scr{P}(S/R)$ to be the submonoid of  $\Scr{S}(S/R)$, consisting entirely of classes $\objeto{[P]}{\phi}{[X]}$, where $[P] \in \Picar{R}$ and $[X] \in \Picar{S}$. By Lemma \ref{lema:31}, the class of an object $\objeto{P}{\phi}{X}$ belongs to $\Scr{P}(S/R)$ provided only $\bara{\phi}_l$ (or $\bara{\phi}_r$) is an isomorphism, see equation \eqref{Phi-l-r}.
\end{definition}

\begin{proposition}\label{pro:32}
Let $R \subseteq S$ be an extension of rings with the same set of local units, and consider $\Scr{P}(S/R)$ as in Definition \ref{PSR}. Then $\Scr{P}(S/R)$ is the units group of the monoid $\Scr{S}(S/R)$.
\end{proposition}
\begin{proof}
It is clear from definitions that if $\Scr{P}(S/R)$ inherits a group structure from $\Scr{S}(S/R)$, then it coincides with the whole units group of $\Scr{S}(S/R)$. Thus we only need to show that $\Scr{P}(S/R)$ is a group with multiplication induced from $\Scr{S}(S/R)$. Since $\Scr{P}(S/R)$ is clearly stable under multiplication and contains the neutral element $\objeto{[R]}{\iota}{[S]}$, it suffices to show that each element $\objeto{[P]}{\phi}{[X]}$ has an inverse in $\Scr{P}(S/R)$. To this end set $P^*=\hom{R}{P}{R}$ and $X^*=\hom{S}{X}{S}$ which we consider canonically as an $R$-bimodule and $S$-bimodule, respectively. Let $[Q]=[P]^{-1}$ in $\Picar{R}$ and $[Y]=[X]^{-1}$ in $\Picar{S}$, we know by equation \eqref{Theta} that there exist natural isomorphisms $\Theta_{-}: -\tensor{R}Q \overset{\cong}{\to} \hom{R}{P}{-}R$ and $\Gamma: -\tensor{S}Y \overset{\cong}{\to} \hom{S}{X}{-}S$. In particular, we have $P^*R \cong Q$ as unital $R$-bimodules and  $X^*S \cong Y$ as unital $S$-bimodules. Define $\objeto{Q}{\psi}{Y}$ an object of $\Scr{M}(S/R)$, where $\psi$ is given by the following composition of $R$-bilinear maps $$\xymatrix{ \psi : Q \ar@{->}^-{\cong}[rr] & & R\tensor{R}Q \ar@{->}^-{\Theta_R}_-{\cong}[rr] &&  P^*R \ar@{->}^-{\underline{\phi^*}}[rr] & & X^*R =X^*S \ar@{->}^-{\Gamma_S^{-1}}_-{\cong}[rr] & & Y}$$ where
$\underline{\phi^*}$ is the restriction of the map $$\xymatrix@R=0pt{  P^* \ar@{->}^-{\phi^*}[rr] & & X^* \\ \sigma \ar@{|->}[rr] & & (\sigma\tensor{R}S) \circ (\bara{\phi}_r){}^{-1}.}$$ In this way, we have a chain of $S$-$R$-bilinear isomorphisms
\begin{footnotesize}
$$\xymatrix@R=27pt@C=50pt{ S\tensor{R}Q \ar@{-->}|-{\overline{\psi_r}}[rrdddd] \ar@{->}^-{\Theta_{S}}[r] & \hom{R}{P}{S}R \ar@{->}^-{\cong}[r] & \hom{R}{P}{\hom{S}{S}{S}R}R \ar@{->}|-{\cong}[d] \\ & & \hom{S}{P\tensor{R}S}{S}R \ar@{->}|-{\cong}[d]  \\ & & \hom{S}{X}{S}R \ar@{=}[d] \\  & & \hom{S}{X}{S}S \ar@{->}|-{\Gamma_{S}^{-1}}[d] \\  & & Y }$$ 
\end{footnotesize}
whose composition coincides exactly with $\overline{\psi_r}$, since we have $\Theta_S(s\tensor{R}q) = s \,\Theta_R(e \tensor{R}q)$ for every element $q \in Q$ and $s \in S$ with unit $e \in R$. We claim that $\objeto{[Q]}{\psi}{[Y]}$ is the inverse of $\objeto{[P]}{\phi}{[X]}$ in $\Scr{P}(S/R)$. For this we only need to show that the following diagram is commutative
$$\xymatrix@C=50pt{ P\tensor{R}Q \ar@{->}_-{\fk{l}}^-{\cong}[d] \ar@{->}^-{\phi\tensor{R}\psi}[r] & X\tensor{R}Y \ar@{->}^-{\omega_{X,Y}}[r] & X\tensor{S}Y \ar@{->}^-{\fk{l}'}_-{\cong}[d] \\ R \ar@{->}|-{[\iota]}[rr] & & S,}$$ where $\fk{l}, \fk{l}'$ are defined as in equation \eqref{Picard-iso}. So let $p \in P$ and $q \in Q$ both with unit $e \in \Sf{E}$, and set
$$\fk{r}^{-1}(e)=\sum q^{(e)}\tensor{R}p^{(e)},\,\, \fk{l}'{}^{-1}(e)= \sum x_{(e)}\tensor{S}y_{(e)},$$ and  for each of those $x_{(e)}$, we set $\bara{\phi_r}^{-1}(x_{(e)})=\sum p_{x_{(e)}}\tensor{R}s_{x_{(e)}}$. Computing the image, we get
\begin{eqnarray*}
  \fk{l}' \circ \omega_{X,Y} \circ (\phi\tensor{R}\psi) (p\tensor{R}q) &=&  \sum \fk{l}'\lr{\phi(p)\tensor{S}e \, \fk{r}\lr{q^{(e)}\fk{l}(p^{(e)}\tensor{R}q)\tensor{R}p_{x_{(e)}}} s_{x_{(e)}} y_{(e)}  } \\
   &=& \sum \fk{l}'\lr{\phi(p)e\, \fk{r}\lr{q^{(e)}\fk{l}(p^{(e)}\tensor{R}q)\tensor{R}p_{x_{(e)}}} s_{x_{(e)}} \tensor{S} y_{(e)} }.
\end{eqnarray*}
On the other hand, we have
\begin{eqnarray*}
\phi(p)e\, \fk{r}\lr{q^{(e)}\fk{l}(p^{(e)}\tensor{R}q)\tensor{R}p_{x_{(e)}}}   &=& \phi\lr{p \,e\, \fk{r}\lr{q^{(e)}\fk{l}(p^{(e)}\tensor{R}q)\tensor{R}p_{x_{(e)}}} } \\
   &=& \phi\lr{p \,\fk{r}\lr{q^{(e)}\fk{l}(p^{(e)}\tensor{R}q)\tensor{R}p_{x_{(e)}}}} \\
   &=& \phi \circ (P\tensor{R}\fk{r}) \circ (P\tensor{R}Q\tensor{R}\fk{l}\tensor{R}P)  \\ &\,\,& \circ  (P\tensor{R}\fk{r}^{-1}\tensor{R}Q\tensor{R}P) (p\tensor{R}e\tensor{R}q\tensor{R}p_{x_{(e)}} )  \\
   &\overset{\eqref{lr-1}}{=}&  \phi \circ (\fk{l}\tensor{R}P) \lr{p\tensor{R}q\tensor{R}p_{x_{(e)}}} \\
   &=& \fk{l}(p\tensor{R}q) \phi(p_{x_{(e)}}),
\end{eqnarray*}
substituting this in the last computation, we get
\begin{eqnarray*}
  \fk{l}' \circ \omega_{X,Y} \circ (\phi\tensor{R}\psi) (p\tensor{R}q) &=&  \sum \fk{l}'\lr{\fk{l}(p\tensor{R}q)\phi(p_{x_{(e)}})s_{x_{(e)}}\tensor{S}y_{(e)}} \\ &=& \sum \fk{l}'\lr{\fk{l}(p\tensor{R}q)\bara{\phi}_r(p_{x_{(e)}}\tensor{R}s_{x_{(e)}})\tensor{S}y_{(e)}} \\ &=& \sum \fk{l}'\lr{\fk{l}(p\tensor{R}q)\bara{\phi}_r\circ \bara{\phi}^{-1}_r(x_{(e)})\tensor{S}y_{(e)}} \\ &=& \sum \fk{l}' \lr{\fk{l}(p\tensor{R}q)x_{(e)}\tensor{S}y_{(e)}} \\ &=& \sum \fk{l}(p\tensor{R}q) \fk{l}'(x_{(e)}\tensor{S}y_{(e)}) \\ &=& \fk{l}(p\tensor{R}q)\, e\,\,=\,\, \fk{l}(p\tensor{R}q).
\end{eqnarray*}
This shows that the above diagram is actually commutative, and this finishes the proof.
\end{proof}

\begin{remark}
In contrast with the unital case, and as we have seen in this Section, the construction of the group $\Scr{P}(S/R)$ is a complicated task. This is due in fact that the inverse of an element $[P] \in \Picar{R}$  (resp. $[X] \in \Picar{S}$) is not represented by the right dual $R$-module $P^*$ (resp. $X^*$). Our approach can be seen also as a complete construction even in the case of extension of rings with unit.
\end{remark}

\section{The third exact sequence of groups}\label{Sect-IV}
Let $R \subseteq S$ be an extension of rings with a same set of local units, and $\Scr{P}(S/R)$ the associated group of Section \ref{Sect-III}. Consider the following maps $$\Scr{D}': \Inv{R}{S} \longrightarrow \Scr{P}(S/R),\quad \lr{X \longmapsto \left( \objeto{[X]}{\subseteq}{[S]} \right) },$$ and
\begin{eqnarray*}
\Scr{O}_r:  \Scr{P}(S/R) \longrightarrow \Picar{S}, & & \lr{(\objeto{[P]}{\phi}{[X]}) \longmapsto [X] },\\
\Scr{O}_l:  \Scr{P}(S/R) \longrightarrow \Picar{R}, & & \lr{(\objeto{[P]}{\phi}{[X]}) \longmapsto [P] }.
\end{eqnarray*}
It is not hard to see that these maps are in fact morphisms of groups. Furthermore, we have

\begin{proposition}\label{pro:IV-1}
Let $R \subseteq S$ be an extension of rings with the same set of local units. Then there is an exact sequence of groups:
$$\xymatrix{1 \ar@{->}[r] & {\rm Ker}(\Scr{D}) \cap {\rm Ker}\widehat{(-)} \ar@{->}[r] & {\rm Ker}\widehat{(-)} \ar@{->}^-{\Scr{D}_{/}}[r] & \Inv{R}{S} \ar@{->}^-{\Scr{D}'}[r] & \Scr{P}(S/R) \ar@{->}^-{\Scr{O}_r}[r] & \Picar{S}, }$$ where $\widehat{(-)}$ is the map of Lemma \ref{lema:22}, and $\Scr{D}_{/}$ is the restriction of the map $\Scr{D}$ defined in Lemma \ref{lema:2}.
\end{proposition}
\begin{proof}
Let $X \in \Inv{R}{S}$ such that $(\objeto{[X]}{\subseteq}{[S]}) = (\objeto{[R]}{\iota}{[S]})$. That is, there exist two isomorphisms $\gamma:X \to R$ and $\lambda: S \to S$,  $R$-bilinear and $S$-bilinear, respectively, rendering commutative the following diagram
$$\xymatrix@C=60pt{ X \ar@{->}_-{\gamma}^-{\cong}[d] \ar@{->}|-{\subseteq}[r] & S  \ar@{->}^-{\lambda}_-{\cong}[d] \\ R \ar@{->}|-{\iota}[r] &  S.}$$ Thus $\gamma(X)=\lambda(X)=R$. So we have $\lambda \in \Aut{S-R}{S}$ with $\lambda^{-1}(R)=X$. Since $\lambda$ is $S$-bilinear, we have $\lambda \in {\rm Ker}\wid{(-)}$. We then get ${\rm Ker}(\Scr{D}') \subseteq {\rm Im}(\Scr{D}_{/})$. Conversely, let $X \in {\rm Im}(\Scr{D}_{/})$, that is, there exists by Proposition \ref{pro:21} $\lambda \in \Aut{S-S}{S}$  such that $X=\lambda^{-1}(R)$. This means that the following diagram
$$\xymatrix@C=60pt{ X \ar@{-->}_-{\lambda}[d] \ar@{->}|-{\subseteq}[r] & S  \ar@{->}^-{\lambda}[d] \\ R \ar@{->}|-{\iota}[r] &  S.}$$ is commutative. Whence, $(\objeto{[X]}{\subseteq}{[S]}) \,=\, (\objeto{[R]}{\iota}{[S]})$, and the converse inclusion is fulfilled. Now, let $(\objeto{[P]}{\phi}{[X]}) \in \Scr{P}(S/R)$ such that $\Scr{O}_r\lr{\objeto{[P]}{\phi}{[X]}}=[X]=[S]$. Thus $X\cong S$ as an $S$-bimodule. Denote by $V$ the copy of $\phi(P)$ in $S$. Since by Lemma \ref{lema:31} $\phi$ is injective, we have $[V]=[P]$ in $\Picar{R}$, and so $V \in \Inv{R}{S}$. Hence $\Scr{D}'(V)= (\objeto{[V]}{\subseteq}{[S]}) = (\objeto{[P]}{\phi}{[X]})$. This shows that ${\rm Im}(\Scr{D}') \subseteq {\rm Ker}(\Scr{O}_r)$. The converse inclusion is easy to check.
\end{proof}

\section{The fourth exact sequence of groups}\label{Sect-V}
Let $R \subseteq S$ be an extension of rings with the same set of local units, and consider the associated groups  $\Aut{R-rings}{S}$ and $\Scr{P}(S/R)$ as were defined in Sections \ref{Sect-II} and \ref{Sect-III}, respectively. Define the map
$$ \mathscr{E}: \Aut{R-rings}{S} \longrightarrow \Scr{P}(S/R), \quad \lr{ \gamma \longmapsto ( \objeto{[R]}{\iota_{\gamma}}{[S_{\gamma}]}) },$$ where $S_{\gamma}$ is the unital $S$-bimodule of Section \ref{Sect-II}. That is, $S_{\gamma}$ is the $S$-bimodule with underlying abelian group $S$ and with biactions: $$ s'\,t_{\gamma}\, s \,\,=\,\, (st\gamma(s))_{\gamma}, \quad \text{ for every } t_{\gamma} \in S_{\gamma}, \text{ and } s,s' \in S.$$ The map $\iota_{\gamma}$ is the canonical inclusion of $R$-bimodules. Using the isomorphisms stated in the proof of Lemma \ref{lema:21}, we can easily check that the map $\mathscr{E}$ is in fact a morphism of groups.

\begin{proposition}\label{pro:41}
Let $R \subseteq S$ be an extension of rings with a same set of local units. Then there is an exact sequence of groups 
$$\xymatrix{ 1 \ar@{->}^-{}[r] & {\rm Ker}(\Scr{D}) \cap {\rm Ker}\wid{(-)} \ar@{->}^-{}[r] & {\rm Ker}(\Scr{D}) \ar@{->}^-{\wid{(-)}_{/}}[r] & \Aut{R-rings}{S} \ar@{->}^-{\Scr{E}}[r] & \Scr{P}(S/R) \ar@{->}^-{\Scr{O}_l}[r] & \Picar{R}, }$$ where $\wid{(-)}_{/}$ is the restriction of the map $\wid{(-)}$ defined in Lemma \ref{lema:22}.
\end{proposition}
\begin{proof}
Let $\gamma  \in \Aut{R-rings}{S}$ such that $(\objeto{[R]}{\iota_{\gamma}}{[S_{\gamma}]}) = (\objeto{[R]}{\iota}{[S]})$ in the group  $\Scr{P}(S/R)$. So we have a commutative diagram
$$\xymatrix@C=60pt{ R \ar@{->}_-{\cong}[d] \ar@{->}|-{\iota}[r] & S  \ar@{->}^-{\lambda}_-{\cong}[d] \\ R \ar@{->}|-{\iota_{\gamma}}[r] &  S_{\gamma},}$$ where $\lambda$ is an isomorphism of $S$-bimodules. The commutativity of this diagram ensures that $\lambda^{-1}(R)=R$. Now, given  $s \in S$  with unit $e$, we have
$$\lambda(es)\,=\, \lambda(e) \gamma(s)\,=\,s\lambda(e), \quad \text{ in the bimodule } S_{\gamma}.$$
Hence $\gamma(s) \,=\, \lambda^{-1}(e) s\lambda(e)$. Therefore, $\gamma \,=\, \wid{\lambda}$ with $\lambda^{-1}(R)=R$. That is, $\gamma \in {\rm Im}(\wid{(-)}_{/})$, and so ${\rm Ker}(\Scr{E}) \subseteq  {\rm Im}(\wid{(-)}_{/})$. Conversely, if $\gamma \in \Aut{R-rings}{S}$ with $\gamma \,=\, \wid{\lambda}$, for some $\lambda \in {\rm Ker}(\Scr{D})$, i.e. $\lambda \in  \Aut{S-R}{S}$ such that $\lambda^{-1}(R)=R$. Then one can easily see that
$$\xymatrix@C=60pt{ R \ar@{->}_-{\lambda}^-{\cong}[d] \ar@{->}|-{\iota}[r] & S  \ar@{->}^-{\gamma}_-{\cong}[d] \\ R \ar@{->}|-{\iota_{\gamma}}[r] &  S_{\gamma},}$$ is a commutative diagram, and thus $(\objeto{[R]}{\iota}{[S]})=(\objeto{[R]}{\iota_{\gamma}}{[S_{\gamma}]})$ in $\Scr{P}(S/R)$. This shows the exactness at the third term of the stated sequence.

Now, let $\objeto{[P]}{\phi}{[X]} \in \Scr{P}(S/R)$ such that $[P]=[R]$. So there is an $R$-bilinear isomorphism $f : R \to P$, which leads to the following $R$-$S$-bilinear and $S$-$R$-bilinear isomorphisms
$$\xymatrix@C=40pt{\alpha: S \ar@{->}^-{\cong}[r] & R\tensor{R}S \ar@{->}^-{f \tensor{R}S}[r] & P \tensor{R}S \ar@{->}^-{\bara{\phi}_r}[r] & X,}$$ $$\xymatrix@C=40pt{\beta: S \ar@{->}^-{\cong}[r] & S\tensor{R}R \ar@{->}^-{S \tensor{R}f}[r] & S \tensor{R}P \ar@{->}^-{\bara{\phi}_l}[r] & X.} $$ If we set $\gamma:= \beta^{-1} \circ \alpha: S \to S$, then by definitions we get
\begin{equation}\label{Ec:gamma}
\gamma(s)\, \phi(f(e))  \,\,=\,\,  \phi(f(e))\, s,\quad \text{ for every } s \in S \text{ with unit } e.
\end{equation}
For every pair of elements $s, t \in S$ with common unit $e$, we obtain
\begin{eqnarray*}
  \alpha(ts) &=& \phi(f(e))\,ts  \\
   &\overset{\eqref{Ec:gamma}}{=}& \gamma(t) \phi(f(e)) s  \\
   &\overset{\eqref{Ec:gamma}}{=}&  \gamma(t) \gamma(s) \phi(f(e)),   \quad e \in Unit\{\gamma(t), \gamma(s)\}, \text{ since } \gamma \text{ is } R\text{-linear} \\ &=&  \bara{\phi}_l (\gamma(t)\gamma(s) \tensor{R} f(e)) \\
   &=& \beta( \gamma(t) \gamma(s)),
\end{eqnarray*}
whence $\gamma(ts)\,=\, \gamma(t) \gamma(s)$. This implies that $\gamma \in \Aut{R-rings}{S}$, since $\gamma$ is $R$-bilinear. Using again equation \eqref{Ec:gamma}, we can easily check that the map $\beta: S_{\gamma} \to S$ is in fact an $S$-bilinear isomorphism. Furthermore, the diagram
$$\xymatrix@C=60pt{ R \ar@{->}_-{f}^-{\cong}[d] \ar@{->}|-{\iota_{\gamma}}[r] & S_{\gamma}  \ar@{->}^-{\beta}_-{\cong}[d] \\ P \ar@{->}|-{\phi}[r] &  X,}$$ is commutative. This shows that $\Scr{E}(\gamma)=(\objeto{[P]}{\phi}{[X]})$. Therefore, ${\rm Ker}(\Scr{O}_l) \subseteq {\rm Im}(\Scr{E})$. The reciprocal inclusion is clear. This completes the exactness of the stated sequence, since the exactness in the second term is obvious.
\end{proof}


The following figure displays the four exact sequences of groups associated to the extension $R \subseteq S$ of rings with the same set of local units.

\begin{footnotesize}
\[
\xymatrixrowsep{8pt} \xymatrixcolsep{15pt}
 \xymatrix{
 & & & 1 \ar@{~>}@/_0.6pc/[dr] & & & & & & \\ &  & &  & {\rm Ker}(\Scr{D}) \ar@{~>}@/_2pc/[dddr] \ar@{.>}_-{}@/^2.5pc/[rr] & & \Aut{R-rings}{S} \ar@{.>}_-{\Scr{E}}@/_1.7pc/[dddr] \ar@{-->}^-{[S_{-}]}@/^2.5pc/[drrr] & & & \\ 1 \ar@{->}_-{}@/^1.7pc/[ddrr] & & &  & &  & & & &\Picar{S} \\ & & & & & & & & & \\ & & {\rm Ker}\wid{(-)} \cap {\rm Ker}(\Scr{D})  \ar@{.>}_-{}@/^2pc/[uuurr] \ar@{->}_-{}@/_2pc/[dddrr] &  &  & \Aut{S-R}{S} \ar@{~>}^-{\Scr{D}}@/_1pc/[dddr] \ar@{-->}^-{\wid{(-)}}@/^1pc/[uuur] &  & \Scr{P}(S/R) \ar@{.>}_-{\Scr{O}_l}@/^1pc/[ddrr] \ar@{->}_-{\Scr{O}_r}@/_1pc/[uurr] & & \\ & & & & & & & & & \\  1 \ar@{.>}_-{}@/_1.7pc/[uurr] & & &  & & & & & & \Picar{R} \\  & & &  & {\rm Ker}\wid{(-)} \ar@{-->}@/^2pc/[uuur] \ar@{->}_-{}@/_2.5pc/[rr] & & \Inv{R}{S} \ar@{~>}_-{[-]}@/_2.5pc/[urrr] \ar@{->}_-{\Scr{D}'}@/^1.7pc/[uuur] & & & \\ &  & & 1 \ar@{-->}@/^0.6pc/[ur] & & & & & & \\ }
\]
\end{footnotesize}

\end{document}